\documentclass[12pt]{amsart}

\usepackage{latexsym,amssymb,color}
\newtheorem{theorem}{Theorem}
\newtheorem{corollary}{Corollary}
\newtheorem{lemma}{Lemma}
\newtheorem{defn}{Definition}
\theoremstyle{remark}
\newtheorem{rem}{Remark}
\def\Map{\mathop{\rm Map}}
\def\Hom{\mathop{\rm Hom}}
\def\Div{\mathop{\rm Div}}
\def\div{\mathop{\rm div}}
\def\poldiv{\mathop{\rm div_{\infty}}}

\def\ch{\mathop{\rm char}}
\def\done{\hfill{$\square$}}

\newcommand{\set}[1]{\left\{#1\right\}}

\newcommand{\Comp}{\mathbb C}

\newcommand{\Intg}{\mathbb Z}
\newcommand{\Nat}{\mathbb N}

\newcommand{\res}{\operatorname{Res}}

\newcommand{\eps}{\varepsilon}
\newcommand{\To}{\longrightarrow}
\newcommand{\CP}{{\Comp\mbox{\rm P}}}
\title[Elliptic genus of complete
intersections]{Mirror symmetry formulae for the elliptic genus of
complete intersections}
\author{Vassily Gorbounov and Serge Ochanine}
\address{Department of Mathematics, University of Kentucky,
Lexington, KY 40506-0027}
\date{\today}
\thanks{The first author is partially supported by the NSF}

\begin{document}

\begin{abstract} In this paper we calculate the elliptic genus of
certain complete intersections in products of projective spaces. We
show that it is equal to the elliptic genus of the Landau-Ginzburg
models that are, according to Hori and Vafa, mirror partners of
these complete intersections. This provides additional evidence of
the validity of their construction. \end{abstract}

\maketitle

Understanding the mathematics behind Quantum and in particular
Conformal Field Theories has been a challenge for more then twenty
five years. The usual ground for mathematical interpretations of
Quantum Field Theory predictions has been the Topological Quantum
Field Theory, which is a certain reduction of the genuine Quantum
Field Theory. There have been significant mathematical advances in
this area. Mirror symmetry is among the major motivations behind
these advances. Vaguely stated, mirror symmetry is a duality between
complex and symplectic geometry. As was originally discovered by
physicists \cite{Cand}, among its specific manifestations is a
striking connection between the``number of curves of a given genus''
in a symplectic manifold and the periods of a holomorphic form of a
different complex manifold.  Ever since, the``explicit''
construction of the mirror partner for a given manifold has become
the central task for the mathematicians. There is a vast amount of
work done in this direction.  In this paper, we  follow the line
pioneered by Gepner \cite{G} and developed  by Vafa \cite{V} who
discovered that the Conformal Field Theory defined by a manifold,
the so-called the sigma model, can be identical to the Conformal
Field Theory of an a different type, the so-called Landau-Ginsburg
model. This point of view was further developed by Witten \cite{W1}.
As far as the application of this idea to mirror symmetry is
concerned, the major reference for the purposes of this paper is the
work of Hori and Vafa \cite{HV}. They showed that the mirror partner
of a large class of manifolds turns out to be a Landau-Ginsburg
model of some kind, or its orbifold. The``proof'' that these models
form a mirror pair would consist of picking an invariant which is
known to be identical for mirror partners and by a calculation
showing that it is indeed the same for given hypothetical mirror
partners. In \cite{HV} such an invariant is given by the periods of
a holomorphic form on a manifold and the so called BPS masses in the
Landau-Ginsburg model. Of course the identities implied by mirror
symmetry in Topological Quantum Field Theory are a reduction of
stronger identities in the original Quantum Field Theory. The
rigorous mathematical structure behind Quantum Field Theory is not
known at the present time, so exploring the mathematical
consequences of mirror symmetry at this higher level is difficult.
Some time ago, Malikov, Schechtman and Vaintrob \cite{MSV}
introduced a construction of a mathematical approximation to the
structure of Quantum Field Theory defined by a manifold. It is so
called the Chiral de Rham complex. Mathematically, these ideas were
further developed in \cite{GMS}, \cite{GMSIII},  \cite{B},
\cite{BL}, \cite{KV}. The relevance of this construction to physics
was explained only recently in \cite{W2}, \cite{Kap}. Despite being
only an approximation to the genuine Quantum Field Theory, the
Chiral de Rham complex carries features not available in the
Topological Quantum Field Theory. One of them is the elliptic genus.
Introduced in mathematics by in \cite{O}, it has been almost
immediately connected to quantum physics \cite{W3},\cite{W5} and
later shown to be identical for mirror partners, thus providing
another test for a pair of manifolds to be mirror partners.
Moreover, as explained in \cite{W4}, there is a physics counterpart
of the elliptic genus in a large class of Conformal Field Theories,
in particular in Landau-Ginsburg models and their orbifolds, and
this physics elliptic genus coincides with the one defined in
topology for the sigma model. Taking on Vafa's and Witten's ideas, a
number of physicists came up with a new type of formulas for the
elliptic genus of some classes of manifolds in terms of their mirror
Landau-Ginsburg partner \cite{KYY}, \cite{NH1}, \cite{NH2},
\cite{GCYLGC} . For about ten years, no mathematical proof of these
formulas was produced. The first paper where such a proof was given
in the case of Calabi-Yau hypersurfaces was \cite{GM}. The scope of
this paper was much broader and the result for the elliptic genus
fell out as a simple consequence of a much deeper connection found
between the Landau-Ginsburg model and the Chiral de Rham complex of
a hypersurface.

The purpose of the present paper is to prove, by more or less
elementary means, that the Landau-Ginsburg mirror partners found in
\cite{HV} for complete intersection in products of projective spaces
have the physics elliptic genus identical to the topological
elliptic genus. It is interesting to note that our result provides
an extra test for the constructions in \cite{HV}, and as such
refines the conditions for the existence of a mirror Landau-Ginsburg
theory as a mirror partner for some complete intersections given in
\cite{HV}. We are planning to return to the case of complete
intersections in more general toric varieties considered in
\cite{HV} in a future work.

Although inspired by physics, with the exception of the introductory
section 0, this paper is mathematically self-contained. It starts
with a brief introduction to Jacobi functions and the correspondent
elliptic genera.  We then proceed to derive the main formula
(theorem 9) using the residue theorem for functions of several
variables.  The formula is derived for elliptic genera of any level
which includes the case of Hirzebruch's level $N$ genus \cite{EGLN},
\cite{MMF}. In section 10, the formula is specialized to the level 2
genus discussed in physics literature. A different proof of the
level 2 case is given in section 11.

The authors benefited from a number of stimulating discussions with
F. Hirzebruch, K. Hori, F. Malikov, D. Zagier. Part of the work was
completed when the first author was visiting the Max Plank Institut
f\"ur Mathematik which provided an excellent stimulating research
environment. The authors are especially grateful to A. Gerasimov for
bringing to their attention the topological aspects of the
Landau-Ginzburg/Calabi-Yau correspondence.

\addtocounter{section}{-1}

\section{Landau Ginsburg orbifolds and their elliptic genus}

In this section we briefly state the results about the elliptic
genus of LG theories and their orbifolds without going into details.
In physics literature the elliptic genus is defined as a character
of an action of an infinite dimensional Lie algebra, namely the so
called $N=2$ algebra. For the purpose of this note, Landau-Ginzburg
field theories are described mathematically by a non-compact
manifold and a function $W$ on it, called \emph{superpotential},
which has isolated singularities. These data are sufficient for
defining and calculating explicitly the elliptic genus of some
relevant for this note LG theories \cite{W4}. An important class of
quantum field theories related to LG is defined by orbifolds of LG
with respect to a finite group action. The elliptic genus of the LG
orbifolds relevant for us was calculated in a number of papers
\cite{NH1}, \cite{NH2}, \cite{GCYLGC}, \cite{KYY} (see also
\cite{GM} for more mathematical approach). Such a genus is either a
Jacobi form or a modular function depending on the data defining the
LG model. In this note we consider LG defined by $(\Comp^N, W)$,
where the superpotential $W$ is a holomorphic and quasi-homogeneous
function of $z_1, z_2, \ldots , z_N$ and its orbifolds with respect
to the action of some finite group of symmetries. The conditions on
the $W$ are the following: It should be possible to assign some
weights $k_i\in \Intg$ to the variables $z_i$ for $i\in \Nat$ and a
degree of homogeneity $d\in \Nat$ to $W$ so that
$$
W(\lambda^{k_i}z_i) = \lambda^d W(z_i)
$$
for all $i$ and $\lambda\in\Comp$. Let $q_i = k_i/d$.

Suppose the potential is invariant under a finite abelian group of
symmetries $G$. Denote by $R_i$ the function on $G$ satisfying
$g(z_i) = \exp(2\pi iR_i(g))z_i$. The invariance of $W$ means of
course that for all $i$ and $g\in G$
$$W(g(z_i)) = W(z_i).$$

The elliptic genus $Ell(q,y)$ of such an LG orbifold defined by the
data $(W,G)$ following \cite{NH1}, \cite{NH2}, is given as follows:

$$
Ell(q,y) = \frac{1}{|G|}\sum_{g_1,g_2\in G} \,\,\prod_{i=1}^{N}
y^{-R_i(g_1)}\frac{\Theta_1((1-q_i)z + R_i(g_1) +
R_i(g_2)\tau|\tau)}{\Theta_1(q_iz + R_i(g_1) + R_i(g_2)\tau |\tau)},
$$
where $y=\exp(2\pi iz)$ and $q=\exp(2\pi i\tau)$.

\emph{Example}. Consider the following LG data taken from
\cite{CVSLGGC}: In the space $\Comp^{9}$ introduce coordinates $X_a$
($a=0,1,2,3$), $Y_b$ ($b=0,1,2$), $Z_c$ ($c=0,1$), and define the
potential by the formula
$$
W = \sum _{r=0}^{1}(X_r^3 + X_rY_r^2 + Y_rZ_r^2) + X_2^3 + X_2Y_2^2
+ X_3^3.
$$

In this case $d=3$ and $q_i=1/3$. It is easy to check that $W$ is
invariant under the following action of $\Intg/12$ with generator
$\omega$:
$$
\omega (X_a, Y_b, Z_c) = (\omega^4 X_a, \omega^{-2} Y_b, \omega
Z_c).
$$
Calculating the appropriate number $R_i(\omega)$ which
define such an action we obtain:
$$
R_a(\omega)=1/3, \quad R_b(\omega)=-1/6,\quad R_c(\omega)=1/12.
$$
Therefore the above formula for the elliptic genus of the LG
orbifold defined by the data $(W, \Intg/12)$ becomes:
\begin{multline*}
Ell(q,y) =\\
\frac{1}{12}\sum_{i,j = 0}^{12}
\bigg(y^{-\frac{j}{3}}\frac{\Theta_1(\frac{2z}{3}-\frac{i+j\tau}{3}
|\tau)}{\Theta_1(\frac{z}{3} + \frac{i+j\tau}{3} |\tau)}\bigg )^4
\bigg (y^{\frac{j}{6}}\frac{\Theta_1(\frac{2z}{3} +
\frac{i+j\tau}{6}|\tau)}{\Theta_1(\frac{z}{3} - \frac{i+j\tau}{6}
|\tau)}\bigg )^3 \bigg (y^{-\frac{j}{12}}\frac{\Theta_1(\frac{2z}{3}
- \frac{i+j\tau}{12} |\tau)}{\Theta_1(\frac{z}{3} +
\frac{i+j\tau}{12}|\tau)}\bigg)^2.
\end{multline*}

It is easy to check that this formula agrees with our formula for
the elliptic genus of the complete intersection in $\CP^3\times
\CP^2\times \CP^1$ given by equations:
$$
\sum_{a=0}^{3} X_a^3 = 0,\quad \sum_{b=0}^{2} X_bY_b^2 = 0,\quad
\sum_{c=0}^{1} Y_cZ_c^2 = 0
$$
(see section 10).

\vspace{12pt}

\begin{quote}
{\em In sections 1--7, $k$ will be an algebraically closed field  of
characteristic $p\geq 0$, $G$ a finite group of order $n$, and
$p\nmid n$.}
\end{quote}

\section{Regular Representations} Consider the $k$-vector space
$V=\Map(G,k)$ of all $k$-valued functions on $G$. Clearly,
$\dim_kV=n$. The group $G$ acts on $V$ by the formula
$$
(g\cdot f)(h) = f(hg), \qquad (g,h\in G, f\in V).
$$
With this action, $V$ is called the {\em regular representation
of} $G$.

We now consider the case when $G$ is abelian and write the group
operation in $G$ additively. Let $\widehat{G} = \Hom(G, k^\ast)$
be the character group of $G$. The following is well known
(cf.~\cite{LRFG}, 2.4):

\begin{theorem}\begin{enumerate}\item $\widehat{G}\subset V$ is a
basis of $V$. \item The one-dimensional subspace $V_\chi$ of $V$
generated by $\chi\in \widehat{G}$ consists of those $f\in V$
which satisfy
$$
f(u+g) \equiv \chi(g)f(u)\qquad\mbox{for all}\quad g\in G.
$$
\item The $G$-modules $V_\chi$ are pairwise non-isomorphic. \item
Every irreducible $G$-module is of degree 1 and is isomorphic to
one of the $V_\chi$.
\end{enumerate} \end{theorem}

\section{Generalized Jacobi Functions} Let now $E$ be an
elliptic curve over $k$, $k(E)$ be the field of rational functions
on $E$, and $G\subset E$ be any finite subgroup of order $n$.  As
above, we assume that $p\nmid n$.

Let $\Delta$ be the divisor $$\Delta=\sum_{g\in G}(g),$$ and let
$\mathcal{L}(G)\subset k(E)$ be the associated vector space:
$$
\mathcal{L}(G)= \left\{f\in k(E)\mid \div(f) \succcurlyeq -
\Delta\right\}.
$$
By the Riemann-Roch theorem, $\dim_k \mathcal{L}(G) = n$. The
group $G$ acts on $E$ by translations and leaves $\Delta$
unchanged. Therefore $\mathcal{L}(G)$ is invariant under the
induced $G$-action on $k(E)$. Thus $\mathcal{L}(G)$ is naturally
an $n$-dimensional representation of $G$.

\begin{theorem} As a $G$-representation, $\mathcal{L}(G)$ is
isomorphic to the regular representation of $G$.\end{theorem}

\proof Let $u_0\in E$ be any point such that $nu_0\neq 0$. Define
a $G$-linear map
$$
\phi : \mathcal{L}(G)\longrightarrow V = \Map(G, k)
$$
by
$$
\phi(f)(g) = f(u_0+g), \qquad (f\in \mathcal{L}(G), g\in G).
$$
This is well-defined, since $u_0\not\in G$ and therefore $u_0+g$
is not a pole of $f$. Since $\mathcal{L}(G)$ and $V$ have the same
dimension $n$, we only need to prove that $\phi$ is injective.

Suppose $f\neq 0$ and $\phi(f)=0$. Then $f(u_0+g)=0$ for all $g\in
G$. Thus $f$ has $\geq n$ zeroes. On the other hand, since $f\in
\mathcal{L}(G)$, it has $\leq n$ poles. Since $\deg\div(f) = 0$,
we conclude that $f$ has a simple pole at each $g\in G$ and a
simple zero at each $u_0+g$. I.e.
$$
\div(f) = - \sum_{g\in G}(g) + \sum_{g\in G}(u_0+g).
$$
Thus the image of $\div(f)$ under Abel's map
$$
\Div(E)\longrightarrow E
$$
is $nu_0\neq 0$. The contradiction shows that $f=0$.\done

\begin{corollary} For each $\chi\in \widehat{G}$, there is a
non-zero function $f_\chi \in \mathcal{L}(G)$ satisfying
$$
f_\chi(u+g)\equiv \chi(g)f_\chi(u)
$$
for all  $g\in G$. This function is determined uniquely up to a
non-zero multiplicative constant. \end{corollary}

\begin{defn}  We call $f_\chi$ a (generalized) Jacobi function
belonging to $\chi\in \widehat{G}$.\end{defn}

\begin{rem} The classical Jacobi functions correspond  to the
case where $G$ is a cyclic group of order 2($\ch k \neq 2$).
\end{rem}

\section{The Divisor of a Jacobi Function} We now describe
explicitly the divisor of a Jacobi function.

\begin{theorem}A non-zero function $f\in \mathcal{L}(G)$
is a Jacobi function if and only if $\div(f)$ is invariant under
translations by elements of $G$.\end{theorem}

\proof If $f$ is a Jacobi function, the formula
$$
f_\chi(u+g)\equiv \chi(g)f_\chi(u),
$$
shows that the functions
$$
u\longmapsto f(u) \quad \mbox{and}\quad u\longmapsto f(u+g)
$$
have the same divisor. Thus $\div(f)$ is invariant by $G$.

Conversely, let $f$ be a non-constant function in
$\mathcal{L}(G)$, and suppose $\div(f)$ is invariant by
translation. Then there is a constant $\chi(g)\in k^\ast$ such
that
$$
f_\chi(u+g)\equiv \chi(g)f_\chi(u),
$$
for all $u\in E$. Since
$$
\chi(g_1+g_2)f(u)=f(u+g_1+g_2)=\chi(g_2)f(u+g_1)=\chi(g_1)
\chi(g_2)f(u),
$$
we have $\chi(g_1+g_2)=\chi(g_1)\chi(g_2)$. Thus $\chi\in
\widehat{G}$, and $f$ is a Jacobi function belonging to
$\chi$.\done

It is now easy to describe the divisors of the Jacobi functions.
Let $E_n$ be the group of $n$-division points on $E$. Since $n$ is
the order of $G$, we have $G\subset E_n$. For every coset
$\gamma\in E_n/G$, let
$$
\Gamma_\gamma = \sum_{r\in \gamma}(r) \in \Div(E).
$$

\begin{theorem} A non-constant function $f\in \mathcal{L}(G)$
is a Jacobi function if and only if
$$
\div(f) = -\Delta +\Gamma_\gamma
$$
for some non-zero $\gamma\in E_n/G$.\end{theorem}

\proof Using theorem 3, we only need to describe the non-zero
principal divisors $D$ satisfying $D \succcurlyeq -\Delta$ and
invariant under $G$.  Since the polar part of $D$ is non-zero and
invariant under $G$, it must be equal to $-\Delta$. Since
$D+\Delta$ is invariant, and is of degree $n$, it is necessarily
of the form
$$
\sum_{g\in G}(r+g)
$$
for some $r\not\in G$, i.e. of the form $\Gamma_\gamma$.\done

\section{Modulus and Conjugation} Let $f$ be a Jacobi function
with character $\chi$ and divisor $ \div(f) = -\Delta
+\Gamma_\gamma$, $\gamma\in E_n/G$. For each $r\in \gamma$, the
involution $u\mapsto r-u$ takes $\div(f)$ to $-\div(f)$.
Therefore, $f(u)f(r-u)$ is a non-zero constant that we call the
\emph{modulus} of $f$ and designate $c(r)$.  Thus $c:\gamma\To
k^\times$.  It is easy to check that
$$
c(r+g) = \chi(g)c(r)
$$
for all $g\in G$.

Similarly, the involution $u\mapsto -u$ takes $-\Delta
+\Gamma_\gamma$ to $-\Delta +\Gamma_\gamma$.  This $u\mapsto
f(-u)$ is also a Jacobi function with character $\chi^{-1}$. In
particular, when $n=2$, $E_n/G$ has order 2, and $\chi^{-1} =
\chi$. It follows that $f(-u) = af(u)$ where $a^2=1$. Thus $f$ is
either even or odd. Since $f$ has a pole of order 1 at $O\in E$,
it has to be odd.

\section{Jacobi Functions on the Tate Curve} Consider a local
field $k$ complete with respect to a discrete valuation $v$, and
let $q\in k^\times$ be any element satisfying $v(q)<1$.  It is
well known (cf. \cite{AEC}, Appendix C, \S 14) that
$E=k^\times/q^\Intg$ can be identified with the Tate curve
$$
y^2+xy = x^3 +a_4x+a_6,
$$
where
$$
a_4 = \sum_{m\geq 1} (-5m^3)\frac{q^m}{1-q^m}
$$
and
$$
a_6 = \sum_{m\geq 1}
\left(-\frac{5m^3+7m^5}{12}\right)\frac{q^m}{1-q^m}.
$$

Let $G\subset E$ be a finite subgroup of order $n$.  As usual, we
assume that $\ch k \nmid n$.  Let $\gamma\in E_n/G$ be any
non-trivial coset.  We will construct an explicit Jacobi function
for $G$, whose zeroes are the points of $\gamma$.

If $\alpha\in k^\times$, we will write $\bar{\alpha}$ for its
image in $E=k^\times/q^\Intg$. Choose $g_1=1, g_2, \ldots, g_n\in
k^\times$ so that
$$
G=\set{\bar{g}_1, \bar{g}_2, \ldots, \bar{g}_n}.
$$
Then choose $r_1, r_2, \ldots, r_n\in k^\times$ so that
$$
\gamma=\set{\bar{r}_1, \bar{r}_2, \ldots, \bar{r}_n}
$$
and
$$
\prod_i r_i = \prod_i g_i.
$$
This can be done as follows.  First, choose $r\in k^\times$ so
that $\bar{r}\in E$ represents $\gamma$. Then clearly,
$$
\gamma = \set{\bar{r}\bar{g}_1, \bar{r}\bar{g}_2, \ldots,
\bar{r}\bar{g}_n}.
$$
Since $\bar{r}^n=1$ in $E$, we have $r^n=q^s$ for some
$s\in\Intg$. Let
$$
r_1=q^{-s}r, \qquad r_i = rg_i\quad (i\geq 2).
$$
Then
$$
\prod_i r_i = q^{-s}r^n \left(\prod_i g_i\right) =\prod_i g_i.
$$

Now consider the basic Theta function
\begin{eqnarray*}
\Theta(u) &=& (1-u^{-1})\prod_{k\geq 1} (1-q^ku)(1-q^ku^{-1})\\
&=& \prod_{k\geq 1} (1-q^ku)\prod_{k\leq 0}(1-q^{-k}u^{-1}).
\end{eqnarray*}
This is an ``analytic'' function on $k^\times$ that has simple
zeroes at the points of $q^\Intg$ (cf. \cite{ATEFLF}) and
satisfies
$$
\Theta(q^{-1}u)= -u\Theta(u).
$$
Along with $\Theta(u)$ we will consider its translates
$\Theta_\alpha(u)$ ($\alpha\in k^\times$) defined by
$$
\Theta_\alpha(u) = \Theta(\alpha^{-1}u).
$$
The function $\Theta_\alpha(u)$ has simple zeroes at the points of
$\alpha q^\Intg$ and satisfies
$$
\Theta_\alpha(q^{-1}u)= -\alpha^{-1}u\Theta_\alpha(u).
$$
Consider
$$
f(u) = \prod_{i=1}^n \frac{\Theta_{r_i}(u)}{\Theta_{g_i}(u)}.
$$
\begin{theorem} The function $f$ is $q$-periodic and defines a
Jacobi function for $G$ having $\gamma$ as the coset of zeroes.
\end{theorem}

\proof It is enough to prove the first statement since this would
imply that $f$ defines a function on $E$ with simple poles at the
points of $G$ and simple zeroes at the points of $\gamma$.  We
have
$$
f(q^{-1}u) = \prod_{i=1}^n
\frac{\Theta_{r_i}(q^{-1}u)}{\Theta_{g_i}(q^{-1}u)} =
\prod_{i=1}^n
\frac{(-r_iu)\Theta_{r_i}(u)}{(-g_iu)\Theta_{g_i}(u)} = f(u).
$$
\done

\section{A Special Case} We now assume that $q$ is not a root of
1, and choose an arbitrary $n$-th root $q_0$ of $q$. For $n\geq 1$
and $a\in \Intg$, define
$$
\vartheta^n_a(u)=\prod_{\ell\geq 1,\, \ell\equiv a (n)} (1-q^\ell
u)\prod_{\ell\leq 0,\, \ell\equiv a (n)}(1-q^{-\ell}u^{-1}).
$$
With the notations of the previous section, we have $\Theta(u) =
\vartheta^1_0(u)$.

\begin{theorem} The function
$$
f(u) = \frac{\vartheta^n_{-1}(u)}{u\,\vartheta^n_{0}(u)}
$$
defines a Jacobi function on $E=k^\times/q^\Intg$ for the group
$G\subset E$ consisting of the images of the $n$-th roots of 1,
and that has simple zeroes at the images of the $n$-th roots of
$q$.\end{theorem}

\proof  Let $\eps_1=1, \eps_2, \ldots, \eps_n$ be the $n$-th roots
of 1 in $k$. Thus
$$
G = \set{\bar{\eps}_1, \bar{\eps}_2,\ldots, \bar{\eps}_n}, \quad
\gamma = \set{\bar{\eps}_1\bar{q}_0,
\bar{\eps}_2\bar{q}_0,\ldots,\bar{\eps}_n\bar{q}_0}.
$$
Notice that
$$
\prod (\eps_iq_0) = \left(\prod \eps_i\right)q.
$$
Thus we can take
$$
r_1 = q^{-1}q_0, \qquad r_i =\eps_iq_0\quad (i\geq 2).
$$
According to the previous theorem, the function
$$
\prod_{i=1}^n \frac{\Theta_{r_i}(u)}{\Theta_{\eps_i}(u)}
$$
has all the required properties.  It remains to identify the
numerator and denominator explicitly. We will be using the
identities
$$
\prod (1-\eps_it) = \prod (1-\eps^{-1}_it) = 1-t^n,
$$
which are easily proven by noticing that the three polynomials in
$t$ have same degree $n$, same roots, and same constant term 1.

We have:
\begin{eqnarray*}
\prod_i\Theta_{\eps_i}(u) &=& \prod_i (1-\eps_iu^{-1})
\prod_{k\geq 1} \prod_i (1-q^k\eps_i^{-1}u) \prod_{k\geq 1}
\prod_i(1-q^k\eps_iu^{-1})\\
&=& (1-u^{-n})\prod_{k\geq 1} (1-q^{nk}u^n) \prod_{k\geq 1}
(1-q^{nk}u^{-n})\\
&=& \vartheta^n_0(u).
\end{eqnarray*}

Similarly,
\begin{eqnarray*}
\prod_i\Theta_{\eps_iq_0}(u) &=& \prod_i (1-\eps_iq_0u^{-1})
\prod_{k\geq 1} \prod_i (1-q^k\eps_i^{-1}q_0^{-1}u) \prod_{k\geq
1} \prod_i(1-q^k\eps_iq_0u^{-1})\\
&=& (1-qu^{-n})\prod_{k\geq 1} (1-q^{nk-1}u^n) \prod_{k\geq 1}
(1-q^{nk+1}u^{-n})\\
&=& \vartheta^n_{-1}(u).
\end{eqnarray*}

Finally,
$$
\prod_i \Theta_{r_i}(u) =\prod_i\Theta_{\eps_iq_0}(u)
\frac{\Theta_{q^{-1}q_0}(u)}{\Theta_{q_0}(u)}
$$
and
$$
\Theta_{q_0}(u) = \Theta(q_0^{-1}u) = \Theta(q^{-1}(qq_0^{-1}u)) =
-qq_0^{-1}u\,\Theta_{q^{-1}q_0}(u),
$$
i.e.
$$
\frac{\vartheta^n_{-1}(u)}{u\,\vartheta^n_{0}(u)}
$$
is a constant multiple of
$$
\prod_{i=1}^n \frac{\Theta_{r_i}(u)}{\Theta_{\eps_i}(u)}
$$
and therefore has all the required properties. \done

\begin{rem} In the case $n=2$, the previous theorem gives
\begin{eqnarray*}
f(u) &=& \frac{1}{u(1-u^{-2})}\prod_{k\geq 1}
\frac{(1-q^{2k-1}u^2)(1-
q^{2k-1}u^{-2})}{(1-q^{2k}u^2)(1-q^{2k}u^{-2})}\\
&=& \frac{1}{u-u^{-1}}\prod_{k\geq 1} \frac{(1-q^{2k-1}u^2)(1-
q^{2k-1}u^{-2})}{(1-q^{2k}u^2)(1-q^{2k}u^{-2})}. \end{eqnarray*}

The formal substitution $u=e^{z/2}$ leads to
$$
zf(e^{z/2}) =\frac{z/2}{\sinh(z/2)}\prod_{k\geq 1}
\frac{(1-q^{2k-1}e^z)(1-
q^{2k-1}e^{-z})}{(1-q^{2k}e^z)(1-q^{2k}e^{-z})},
$$
which is a familiar expression for the generating function of the
level 2 elliptic genus (see below).
\end{rem}

\section{Modulus and Normalization on the Tate curve} Continuing
with the situation of the preceding section, we are first going to
compute the modulus for the Jacobi function defined by
$$
f(u) = \frac{\vartheta^n_{-1}(u)}{u\,\vartheta^n_{0}(u)},
$$
i.e. compute the constant value of $f(u)f(ru^{-1})$ for $r\in
\gamma$.  We start with $r=q_0$.  The following formulas can be
easily obtained from the definition of $\vartheta^n_a$:
$$
\vartheta^n_0(u^{-1})= -u^n\vartheta^n_0(u), \qquad
\vartheta^n_a(u^{-1})= -u^n\vartheta^n_{-a}(u)\quad (a\not\equiv 0
(n))
$$
and
$$
\vartheta^n_0(q_0u)= -q^{-1}u^{-n}\vartheta^n_1(u), \qquad
\vartheta^n_{-1}(q_0u)=\vartheta^n_0(u).
$$
It follows that
$$
\vartheta^n_0(q_0u^{-1})=-q^{-1}u^n\vartheta^n_{-1}(u), \qquad
\vartheta^n_{-1}(q_0u^{-1}) = -u^n\vartheta^n_0(u).
$$
and therefore
$$
f(q_0u^{-1}) = \frac{qu\vartheta^n_0(u)}{q_0\vartheta^n_{-1}(u)} =
q_0^{n-1}f(u)^{-1}.
$$
Thus
$$
c(q_0) = q_0^{n-1}.
$$
Notice now that the formula for $f$ does not depend on the choice
of $q_0$.  Therefore we have the following
\begin{theorem} For the Jacobi function defined by
$$
f(u) = \frac{\vartheta^n_{-1}(u)}{u\,\vartheta^n_{0}(u)},
$$
and any $n$-th root $r$ of $q$ in $k^\times$, we have
$$
c(r) = r^{n-1}.
$$
\end{theorem}

We now normalize $f$ by requiring that
$$
\res_{u=1}\left( f(u)\frac{du}{u}\right) = \frac{1}{n}.
$$
This normalization is formally equivalent to the requirement that
$zf(e^{z/n})=1+o(z)$.  Since
\begin{eqnarray*}
\frac{(u-1)\vartheta^n_{-1}(u)}{u\,\vartheta^n_{0}(u)} &=&
\frac{u-1}{u(1-u^{-n})}\cdot
\frac{\vartheta^n_{-1}(u)}{\prod_{k\geq
1} (1-q^{nk}u^n)(1-q^{nk}u^{-n})}\\
&=& \frac{u^{n-1}}{1+u+u^2+\cdots +u^{n-1}}\cdot
\frac{\vartheta^n_{-1}(u)}{\prod_{k\geq 1}
(1-q^{nk}u^n)(1-q^{nk}u^{-n})},
\end{eqnarray*}
we see that the normalized Jacobi function is $f(u)N$, where
$$
N =\frac{\prod_{k\geq 1} (1-q^{nk})^2}{\vartheta^n_{-1}(1)}.
$$
The modulus of the normalized function is
$$
c(r)=r^{n-1}N^2.
$$
and satisfies
$$
c(r)^n = q^{n-1}N^{2n}.
$$
In particular, if $n=2$, we have
$$
N= \prod_{k\geq 1} \left(\frac{1-q^{2k}}{1-q^{2k-1}}\right)^2,
$$
and
$$
c(r)^2 = q\prod_{k\geq 1}
\left(\frac{1-q^{2k}}{1-q^{2k-1}}\right)^8,
$$
which is the familiar expression for the modular form $\eps$ for
one of the three level two elliptic genera.

\section{Elliptic Genus} From now on, $k=\Comp$, $E=\Comp/L$
for some lattice $L$, $G=L_0/L$ and $n=[L_0:L]$.  Let $f$ be a
Jacobi function with character $\chi$.  Then $f$ is an elliptic
function with period lattice $L$.  We \emph{normalize} $f$ so that
it has residue $\res_{\,\,z=0}(f)=1$.  Then the Taylor expansion
of $zf(z)$ is a formal power series with constant term 1 and
defines, via the Hirzebruch formalism, a multiplicative genus
$$
\varphi: \Omega^U_*\To \Comp,
$$
which we refer to as the \emph{level $n$ elliptic genus} defined
by $f$. The case $n=2$ is best known. In this case, $zf(z)$ is an
even series in $z$ and $\varphi$ factors through a genus
$$
\varphi: \Omega^{SO}_*\To \Comp.
$$

\section{Complete Intersections} In this section, we will compute
the elliptic genus of complete intersections satisfying a
non-degeneracy condition as a summation over some division points
of $E\times E\times\cdots\times E$.

Let $M=(m_{ij})$ be a $l\times t$ matrix over $\Intg$, and let $P$
be the product
$$
P = \CP^{N_1-1}\times \CP^{N_2-1} \times\cdots\times \CP^{N_t-1}.
$$
For $1\leq j\leq t$, let $\eta_j$ be the pull-back over $P$ of the
canonical line bundle of $\CP^{N_j-1}$. Then for $1\leq i\leq l$,
let $\xi_i$ be the line bundle
$$
\xi_i = \eta_1^{m_{i1}}\otimes \eta_2^{m_{i2}}\otimes \cdots
\otimes \eta_t^{m_{it}}
$$
over $P$. We will write $H_i$ for the stably almost complex
manifold (hypersurface) dual to $\xi_i$ and let $X(M)$ be the
transverse intersection
$$
X(M) = H_1\cap H_2\cap\cdots\cap H_l.
$$

Writing $Mz=(\mu_1z, \mu_2z,\ldots,\mu_lz)$, we define $l$ linear
forms $\mu_i: \Comp^t\To \Comp$, i.e.
$$
\mu_iz = m_{i1}z_1+m_{i2}z_2+\cdots +m_{it}z_t.
$$
Let $\varphi$ be the elliptic genus of level $n$ defined by a
Jacobi function for $G\subset E$ with character $\chi$. The
standard computation using Hirzebruch's formalism leads to the
following

\begin{theorem} $\varphi(X(M))$ is the coefficient of
$z_1^{-1}z_2^{-1}\cdots z_t^{-1}$ in the Laurent expansion at
$z=(0,0,\ldots, 0)$ of
$$
F(z)=F(z_1,z_2,\ldots,z_t)=\frac{f(z_1)^{N_1}f(z_2)^{N_2}\cdots
f(z_t)^{N_t}}{f(\mu_1z)f(\mu_2z)\cdots f(\mu_lz)}.
$$
\done
\end{theorem}

In accordance with the notation $[m]$ for the
multiplication-by-$m$ map on $E$, we will write
$$
[M]: E^t\To E^s, \qquad [\mu_i]: E^t\To E, \qquad [\mu]: E^t\To E
$$
for the maps induced by left multiplication by $M$, by $\mu_i$,
and by $\mu = \sum \mu_i$.

For the rest of this section we will assume that $l=t$ and that
$\det M\neq 0$.  Let $G^t=G\times \cdots \times G\subset E^t$ and
$H\subset E^t$ be the inverse image of $G^t$ under $[M]$. Clearly,
$H$ is a subgroup of $E^t$ containing $G^t$ (since the entries of
$M$ are integers), and it is easy to check that $[H:G^t]=(\det
M)^2$.

For each $j$, $1\leq j\leq t$, choose a zero $r_j\in E$ of $f$ so
that
$$
\div(f) = -\sum_g(g) + \sum_g(r_j+g),
$$
and let $r=(r_1,r_2,\ldots,r_t)\in E^t$.  Also, choose a $s=(s_1,
s_2, \ldots, s_t)\in E^t$ satisfying $[M]s=r$ (this uses $\det
M\neq 0$).

\begin{theorem}  If for every $j$, $\sum_i m_{ij} \equiv N_j \mod
\exp(G)$, we have
$$
\varphi(X(M))=\frac{(-1)^{t+1}}{(\det M)c(r)}\sum_{h\in H/G^t}
\chi(-[\mu]h)f(s+h)^N,
$$
where we use the abbreviations:
$$
c(r) = c(r_1)c(r_2)\cdots c(r_t)
$$
and
$$
f(s+h)^N = f(s_1+h_1)^{N_1}f(s_2+h_2)^{N_2}\cdots
f(s_t+h_t)^{N_t},
$$
and where the sum runs over representatives of the cosets in
$H/G^t$.
\end{theorem}

\begin{rem} The condition $\sum_i m_{ij} \equiv N_j \mod
\exp(G)$ is equivalent to $c_1(X(M))\equiv 0 \mod\exp(G)$.
\end{rem}

The proof of this theorem is based on the following, slightly
modified version of the global residue theorem for functions of
several complex variables as described in chapter 5 of \cite{PAG}
(pp. 655--656).

\begin{theorem} Let $V$ be a compact complex manifold of dimension
$t$, and let $$D_1, D_2, \ldots, D_\nu$$ ($t\leq\nu$) be effective
divisors having the property that the intersection of every $t$ of
them, $D_{n_1}\cap D_{n_2}\cap\cdots\cap D_{n_t}$, is a finite set
of points, whereas the intersection of every $t+1$ of them is
empty. Let $D=D_1 + D_2+\cdots +D_\nu$, and let $\omega$ be a
meromorphic $t$-form on $V$ with polar divisor $D$.  Then for
every $P\in D_{n_1}\cap D_{n_2}\cap\cdots\cap D_{n_t}$ the residue
$\res_P\omega$ is defined and
$$
\sum_P \res_P\omega = 0.
$$
\end{theorem}

The residue $\res_P\omega$ for $P\in D_{n_1}\cap
D_{n_2}\cap\cdots\cap D_{n_t}$ can be defined as follows: choose
local coordinates $(z_1, z_2,\ldots, z_t)$ near $P$ and write
$\omega$ near $P$ as
$$
\omega = \frac{\psi(z)\,dz_1\wedge dz_2\wedge\cdots\wedge dz_t}{
\phi_1(z)\phi_2(z)\cdots\phi_t(z)},
$$
where $\psi, \phi_1, \phi_2,\ldots\phi_t$ are pairwise relatively
prime holomorphic functions. Then
$$
\res_P\omega = \frac{1}{(2\pi i)^t} \int_\Gamma \omega,
$$
where $\Gamma$ is the suitably oriented real cycle
$$
\Gamma = \set{z \mid \phi_1(z)=\varepsilon, \phi_2(z)=\varepsilon,
\ldots, \phi_t(z)=\varepsilon},
$$
for a small $\varepsilon$.

It can be easily verified that if $F(z)$ is a meromorphic function
near $O=(0, 0, \ldots, 0)\in \Comp^t$, then $\res_O
(F(z)\,dz_1\wedge dz_2\wedge\cdots\wedge dz_t$) is the coefficient
of $z_1^{-1}z_2^{-1}\cdots z_t^{-1}$ in the Laurent expansion of
$F(z)$ at $O$.  We will also use the following formula,
generalizing a fact well-known for functions of one variable:

\begin{lemma} If $\omega$ is written near $P$ as
$$
\omega = \frac{\psi(z)\,dz_1\wedge dz_2\wedge\cdots\wedge dz_t}{
\phi_1(z)\phi_2(z)\cdots\phi_t(z)},
$$
and $\det(\partial \phi_i/\partial z_j)(P)\neq 0$, then
$$
\res_P\omega = \frac{\psi(P)}{\det(\partial \phi_i/\partial
z_j)(P)}.
$$
\end{lemma}

We now turn to the proof of theorem 10. In view of theorem 9, and
the above remarks, we need to compute $\res_O \omega$, where
$$
\omega =\frac{f(z_1)^{N_1}f(z_2)^{N_2}\cdots
f(z_t)^{N_t}\,dz_1\wedge dz_2\wedge\cdots\wedge
dz_t}{f(\mu_1z)f(\mu_2z)\cdots f(\mu_lz)}.
$$
We let $V=E^t$ and consider $\omega$ on $V$. According to theorem
11, we have
$$
\sum_P \res_P\omega = 0,
$$
where $P$ runs through the points $P\in G^t$, corresponding to the
poles of $f(z)$ (first kind), and through the points $P\in s+H$
which are the simultaneous zeroes of $f([\mu_1]z), f([\mu_2]z),
\ldots, f([\mu_l]z)$ (second kind).

Notice that the condition
$$
\sum_i m_{ij} \equiv N_j \mod \exp(G)
$$
guarantees that $F(z+g) =F(z)$ for all $g\in G^t$.  It follows
that the contribution of the points of the first kind to the sum
of residues is $n^t\res_O \omega$.

Turning to the points of the second kind, say $P=s+h$, we first
notice that differentiating
$$
f(r_i - z)=\frac{c(r_i)}{f(z)}
$$
with respect to $z$ and taking the limit as $z\rightarrow 0$ gives
$$
f'(r_i)=-c(r_i).
$$
Therefore
$$
\frac{\partial f([\mu_i]z)}{\partial z_j}(s+h) =
f'([\mu_i](s+h))m_{ij} =
f'(r_i+[\mu_i]h))m_{ij}=-\chi([\mu_i]h)c(r_i)m_{ij}.
$$
Thus
$$
\det(\partial f([\mu_i]z)/\partial z_j)(s+h) =
(-1)^tc(r)\chi([\mu]h),
$$
and
$$
\res_{s+h} \omega = \frac{(-1)^t\chi(-[\mu]h)f(s+h)^N}{c(r)}.
$$
Also, notice that this residue remains unchanged when $h$ is
replaced by $h+g$ with $g\in G^t$ (because of the condition $
\sum_j m_{ij} \equiv N_j \mod \exp(G)$).

Theorem 10 is now an immediate consequence of the residue theorem
11.

\section{The Level 2 Case} We now specialize the
above formula for the level 2 elliptic genus. Let $\tau\in
\mathcal{H} = \set{z\in \Comp \mid \Im(z)>0}$,
$L_0=\Intg\oplus\Intg\tau$, $L=\Intg\oplus\Intg(2\tau)$,
$G=\set{0, \tau}$,  and let $\chi: G\To \set{\pm 1}$ be defined by
$\chi(\tau)=-1$. The divisor of the corresponding Jacobi function
$f$ is
$$
\div(f) = -(0)-(\tau)+(1/2)+(1/2+\tau).
$$
By comparing the divisors and the $1/z^2$ terms in the Taylor
expansions at $0$, we easily conclude that
$$
f(z)^2 = \wp(z\vert \tau) - \wp(1/2\,\vert\, \tau)=\wp(z\vert
\tau) - e_1.
$$
We will choose $r=1/2$. Since $f$ satisfies the differential
equation
$$
(f')^2 = f^4 -2\delta f^2+ \eps
$$
(cf.~\cite{EG}), we see that $c(1/2)^2 = f'(1/2)^2 =\eps$, and,
with an appropriate choice of the square root,
$c(1/2)=\sqrt{\eps}$.

Consider first the case of a hypersurface $X(m)\subset \CP^{N-1}$
of degree $m$. Since we can take $s=1/2m$, we have:

\begin{theorem} If $N$ and $m$ have same parity, then
$$
\varphi(X(m)) = \frac{1}{m\sqrt{\eps}} \sum_{0\leq a, b< m} (-1)^a
f^N\left(\frac{1}{2m} + \frac{b}{m} +\frac{a\tau}{m}\right).
$$
\end{theorem}

This is the Eguchi--Jinzenji formula derived in~\cite{GCYLGC}.
Notice that only the case where both $N$ and $m$ are even is of
interest, since, for dimension reasons, the genus vanishes when
$N$ is odd.

Turning now to complete intersections, consider the case of
$X(M)\subset \CP^3\times \CP^2\times \CP^1$, where
$$
M=\begin{pmatrix}
  3 & 0 & 0 \\
  1 & 2 & 0 \\
  0 & 1 & 2
\end{pmatrix}.
$$

If we take $r=(1/2, 1/2, 1/2)$, then one can take $s=(1/6, 1/6,
1/6)$, since
$$
\begin{pmatrix}
  3 & 0 & 0 \\
  1 & 2 & 0 \\
  0 & 1 & 2
\end{pmatrix}\begin{pmatrix}
  1/6 \\
  1/6 \\
  1/6
\end{pmatrix}=
\begin{pmatrix}
  1/2\\
  1/2\\
  1/2
\end{pmatrix}.
$$
Then, noticing that
$$
\begin{pmatrix}
  0 & 1 & 0 \\
  0 & 0 & 1 \\
  1 & -3 & 6
\end{pmatrix}
\begin{pmatrix}
  3 & 0 & 0 \\
  1 & 2 & 0 \\
  0 & 1 & 2
\end{pmatrix}
\begin{pmatrix}
  1 & -2 & 4 \\
  0 & 1 & -2 \\
  0 & 0 & 1
\end{pmatrix}=
\begin{pmatrix}
  1 & 0 & 0 \\
  0 & 1 & 0 \\
  0 & 0 & 12
\end{pmatrix},
$$
we easily conclude that
$$
\frac{1}{12}\begin{pmatrix}
  1 & -2 & 4 \\
  0 & 1 & -2 \\
  0 & 0 & 1
\end{pmatrix}
\begin{pmatrix}
  0 \\
  0 \\
  1
\end{pmatrix}=
\begin{pmatrix}
  1/3 \\
  -1/6 \\
  1/12
\end{pmatrix} = v
$$
generates a cyclic subgroup of order 12 of $H/G^3$, so that one
can take
$$\set{(a+b\tau)v \mid 0\leq a, b< 12}$$ as representatives
of the cosets of $H/G^3$. Also, notice that $[\mu]v = 1$, so
$$\chi(-[\mu](a+b\tau)v)=(-1)^b.$$  We obtain the following
\begin{theorem} For $M$ as above,
$$
\varphi(X(M))=\frac{1}{12\,\eps^{3/2}}\sum_{0\leq a,b<12}(-1)^b
f\left(\frac{1}{6}+\frac{a+b\tau}{3}\right)^4f\left(\frac{1}{6}-\frac{a+
b\tau}{6}\right)^3f\left(\frac{1}{6}+\frac{a+b\tau}{12}\right)^2.
$$
\end{theorem}

The general case can be treated in exactly the same way, since one
can always find invertible over $\Intg$ matrices $A, B$ such that
$AMB$ is diagonal.

\section{A different approach for the level 2 case}  Considering, as
before, a hypersurface  $X(m)\subset \CP^{N-1}$ of degree $m$, and
keeping the assumptions of the previous section, we first note that
$$
\div(f) = -\sum_g (g) + \sum_g (r+g),
$$
so that
$$
\div(f(mz)) = -\sum_h (h) + \sum_h (s+h),
$$
and thus the polar divisor of $1/f(mz)$ is
$$
\poldiv\left(\frac{1}{f(mz)}\right) =-\sum_h (s+h).
$$
To calculate the residue of $1/f(mz)$ at $s+h$, let $t\mapsto
\lambda(t)$ be a small loop going around $0$ in the positive
direction and for $h\in H$, let $\mu(t) =\lambda(t)+s+h$. Then,
using the substitution $u=z+s+h$, we have:
\begin{eqnarray*}
\res_{\,\,z=s+h}\left(\frac{1}{f(mz)}\right)&=&\frac{1}{2\pi i}
\oint_\mu \frac{dz}{f(mz)}\\
&=&\frac{1}{2\pi i}\oint_\lambda \frac{du}{f(mu +r +mh)}\\
&=& \frac{\chi(-mh)}{2\pi i}\oint_\lambda \frac{du}{f(mu +r)}\\
&=& \frac{\chi(-mh)}{2\pi ic(r)}\oint_\lambda f(-mu)du\\ &=&
-\frac{\chi(-mh)}{mc(r)}.
\end{eqnarray*}.

\begin{theorem} We have
$$
\frac{1}{f(mz)} = -\frac{1}{mc(1/2)} \sum_{0\leq a,b <m}
\frac{(-1)^a f'(z)}{f(z) - f(\omega_{a,b})},
$$
where
$$
\omega_{a,b} =\frac{1}{2m} + \frac{b}{m} +\frac{a\tau}{m}.
$$
\end{theorem}

\proof Specializing the above calculations to the level 2 case, we
see that the polar divisor of $1/f(mz)$ is given by
$$
\poldiv\left(\frac{1}{f(mz)}\right) = - \sum_{0\leq a<2m,\, 0\leq
b<m} (\omega_{a,b}),
$$
and that
$$
\res_{\omega_{a,b}}\left(\frac{1}{f(mz)}\right) =
-\frac{(-1)^a}{mc(1/2)}.
$$
Consider the function
$$
h_{a,b}(z) = \frac{f'(z)}{f(z) - f(\omega_{a,b})},
$$
i.e. the logarithmic derivative of $f(z) - f(\omega_{a,b})$. Let
$\sigma$ be the involution $\sigma(z)=\tau -z$.  Then, since $f$
is odd,
$$
f(\tau -z) = -f(-z) = f(z),
$$
i.e. $f\circ \sigma = f$. Also,
$$
\sigma(\omega_{a,b}) = \tau -\frac{1}{2m} - \frac{b}{m}
-\frac{a\tau}{m} = \frac{1}{2m} + \frac{-b-1}{m}
+\frac{(m-a)\tau}{m} = \omega_{m-a, -b-1},
$$
where $m-a$ and $-b-1$ should be interpreted modulo $2m$ and $m$
respectively. As a consequence, we have $h_{a,b} = h_{m-a, -b-1}$.

The function $f(z) - f(\omega_{a,b})$ has two simple poles $0$ and
$\tau$. Thus, it has two zeroes (counted with multiplicities), one
of which is $\omega_{a,b}$. Since $m$ is even, $b\not\equiv -b-1
\mod m$. Thus $\omega_{m-a, -b-1}$ is the other zero of $f(z) -
f(\omega_{a,b})$ and both zeroes are simple. Thus
$$
\poldiv(h_{a,b}) = -(0)-(\tau) -(\omega_{a,b}) - (\omega_{m-a,
-b-1}),
$$
with residue $-1$ at $0$ and $\tau$, and residue $1$ at
$\omega_{a,b}$ and $\omega_{m-a, -b-1}$. Let now
$$
F(z) =  -\frac{1}{mc(1/2)} \sum_{0\leq a<m,\, 0\leq b<m}
\frac{(-1)^a f'(z)}{f(z) - f(\omega_{a,b})} = -\frac{1}{2mc(1/2)}
\sum_{0\leq a<2m,\, 0\leq b<m} (-1)^a h_{a,b}(z).
$$
The possible poles of $F$ are $0$, $\tau$, and $\omega_{a,b}$ with
$0\leq a<2m,\, 0\leq b<m$.  Let's compute the residues of $F$ at
these points. First, we have
$$
\res_0(F) = \res_\tau(F) =  -\frac{1}{2mc(1/2)} \sum_{0\leq
a<2m,\, 0\leq b<m} (-1)^a(-1) = 0,
$$
since $m$ is even.  Then,
$$
\res_{\omega_{a,b}}(F) =-\frac{1}{2mc(1/2)} ((-1)^a+ (-1)^{m-a})=
-\frac{(-1)^a}{mc(1/2)}.
$$
Thus, $1/f(mz)$ and $F(z)$ have the same polar part. It follows
that
$$
\frac{1}{f(mz)} = F(z) +C,
$$
for some constant $C$.  Replacing $z$ with $\sigma(z)=\tau-z$, we
have
$$
\frac{1}{f(m\sigma(z))} = \frac{1}{f(m\tau - mz)}= \frac{1}{f(-
mz)}= -\frac{1}{f(m\sigma(z))},
$$
since $m\tau$ is a period of $f$ for $m$ even. Similarly, since
$f(z)$ is invariant under $\sigma$ and $f'(\sigma(z))=-f'(z)$, we
have $h_{a,b}(\sigma(z))=-h_{a,b}(z)$ and therefore
$F(\sigma(z))=-F(z)$.  It follows that $C=0$.

As a corollary, we obtain a new proof of theorem 11.  Indeed,
according to theorem 8,
$$
\varphi(X(m)) = \frac{1}{2\pi i} \oint_\lambda
\frac{f(z)^N\,dz}{f(mz)},
$$
where $\lambda$ is a small circle around 0 traversed
counter-clockwise.  Thus
$$
\varphi(X(m)) =-\frac{1}{2\pi imc(1/2)}\oint_\lambda \ \sum_{0\leq
a,b <m} \frac{(-1)^a f(z)^Nf'(z)\,dz}{f(z) - f(\omega_{a,b})}.
$$
Changing the variable to $u=f(z)$, we get
$$
\varphi(X(m)) =-\frac{1}{2\pi imc(1/2)}\oint_\Lambda \ \sum_{0\leq
a,b <m} \frac{(-1)^a u^N,du}{u - f(\omega_{a,b})},
$$
where $\Lambda$ is a large circle around 0 traversed
\emph{clockwise}. The Cauchy integral formula gives
$$
\varphi(X(m)) =\frac{1}{mc(1/2)} \sum_{0\leq a,b <m} (-1)^a
f(\omega_{a,b})^N,
$$
which is exactly the formula in theorem 11.

\bibliographystyle{plain}

\end{document}